\documentclass[]{interact}
\usepackage{xcolor}
\usepackage{epstopdf}
\usepackage{natbib}
\bibpunct[, ]{(}{)}{;}{a}{}{,}
\theoremstyle{plain}

\theoremstyle{definition}

\theoremstyle{remark}

\usepackage{subfig}
\usepackage{graphicx}

\usepackage{placeins}
\begin{document}

\articletype{ARTICLE TEMPLATE}

\title{Modelling host population support for combat adversaries}

\author{
\name{
Mathew Zuparic \textsuperscript{a},
Sergiy Shelyag\textsuperscript{b},
Maia Angelova\textsuperscript{b},
Ye Zhu\textsuperscript{b} and 
Alexander Kalloniatis \textsuperscript{a}
}
\affil{
\textsuperscript{a} 
Defence Science and Technology Group, Canberra, ACT 2600 Australia;\\
\textsuperscript{b} 
Data to Intelligence Research Centre, Deakin University, Geelong,
VIC 3220 Australia.
}
}
\maketitle

\begin{abstract}
We consider a model of adversarial dynamics consisting of three populations, labelled Blue, Green and Red, which evolve under a system of first order nonlinear differential equations. Red and Blue populations are adversaries and interact via a set of Lanchester combat laws. Green is divided into three sub-populations: Red supporters, Blue supporters and Neutral. Green support for Red and Blue leads to more combat effectiveness for either side. From Green's perspective, if either Red or Blue exceed a size according to the capacity of the local population to facilitate or tolerate, then support for that side diminishes; the corresponding Green population reverts to the neutral sub-population, who do not contribute to combat effectiveness of either side. The mechanism for supporters deciding if either Blue or Red are too big is given by a logistic-type interaction term. The intent of the model is to examine the role of influence in complex adversarial situations typical in counter-insurgency, where victory requires a genuine balance between maintaining combat effectiveness and support from a third party whose backing is not always assured. 
\end{abstract}



\maketitle


\section{Introduction}
\label{intro}

\begin{figure}[!htb]
\begin{center}
\includegraphics[width=4.cm]{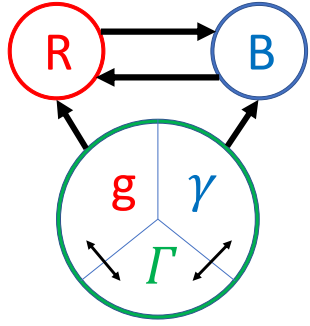}
\caption{Diagram showing \textit{Blue}, \textit{Red} and \textit{Green} populations, and their interactions.}
\label{fig:Model}
\end{center}
\end{figure}

The intent of this work is to explore the mutual influence between \textit{non-combatants} and \textit{combatants} in warfare. To enable this we devise a scenario shown in Figure \ref{fig:Model} which contains three distinct populations, labeled Blue, Red and Green. The Blue and Red forces are homogeneous and in direct opposition with each other. Green represents a non-homogeneous civilian population in proximity to the Blue and Red forces, divided into three sub-populations. Sub-populations $\gamma$ and $g$ support the Blue and Red forces, respectively. The third sub-population $\Gamma$ remains completely neutral, offering no support for either Blue or Red. The heterogeneity of Green is similar to the model proposed in \citep{Atkinson12}, where the authors divided their Blue and Red populations into regions of \textit{supporters} and \textit{opponent sympathisers}. Arguably,
this is a situation often encountered in
counter-insurgency where government forces (in recent history, supported by external actors) contest opponents in a space where members of a population
between the two may be swayed in one direction or the other.

The model dynamics reflects the arrows which bridge the boundaries of the populations in Figure \ref{fig:Model}. Blue and Red affect each other negatively via a \textit{Lanchester}-type interaction. The support offered to Red and Blue from the sub-populations $g$ and $\gamma$ respectively, is determined by the magnitude of the sub-populations. Additionally, the model allows for dynamic movement amongst the three Green sub-populations. The movement from one sub-population to another in this model represents a shift in the Green population's cultural perception of the conflict between Blue and Red, akin to the recent application of the Axelrod cultural model by \citep{Avella14}. The direction and rate of the flow between the Green sub-populations is determined by a \textit{logistic-type} expression; if either Blue or Red become too big the support they derive from their corresponding Green sub-populations begins to wane. 



\subsection{Mathematical modelling of adversarial settings}
\label{mathmod}
An important milestone of adversarial mathematical modelling can be traced back to Verhulst's (\citeyear{Verhulst1838}) formulation of the \textit{logistic equation}:
\begin{equation}
    \dot{P} = \beta P \left(1 - \frac{P}{k} \right), \;\; P(0) = P_0,
    \label{logeq}
\end{equation}
for population $P \in \mathbb{R}$. In Eq.(\ref{logeq}), the tension between the linear growth term with \textit{growth rate} $\beta \in \mathbb{R}_+$, and the quadratic decay term with \textit{carrying capacity} $k \in \mathbb{R}_+$, ensures population equilibrium at the carrying capacity value. The tension between Eq.(\ref{logeq})'s linear growth and nonlinear decay remains a central theoretical mechanism for the modelling of single species \citep{Richards59}, and multi-species dynamics \citep{Royama71}. 


The \textit{force attrition model} of Osipov (\citeyear{Osipov15}) and Lanchester (\citeyear{Lanchester16}) is also a popular mathematical model to understand combat outcomes. The model has two typical variants given by
\begin{eqnarray}
\begin{split}
& \dot{B} = - \beta_{R} R,\;\; \dot{R} = - \beta_B B, \;\; \textrm{square law/direct fire},\\
\textrm{or}\;\; & \dot{B} = - \beta_R BR,\;\; \dot{R}= -\beta_B BR, \;\; \textrm{linear law/indirect fire},
\end{split}
\label{lanch}
\end{eqnarray}
where $B$ and $R$ correspond to the magnitudes of the opposing Blue and Red forces with initial conditions $\{B(0),R(0)\} = \{B_0,R_0\}$, and the parameters $\beta_B$ and $\beta_R$ represent their respective rates of effectiveness. Eq.(\ref{lanch}) was devised as a simple model to understand the interplay between total force numbers and fighting effectiveness in static environments \citep{MacKay06}. Generalisations of the Lanchester equations include Deitchman's \textit{guerrilla warfare} model (\citeyear{Deitchman62}), which mixes the linear and square Lanchester laws to account for the effectiveness of an opponent employing unconventional tactics, and MacKay's \textit{mixed forces} model (\citeyear{MacKay09,MacKay12}) for a general number of $N$ Blue forces vs $M$ Red forces.

Previous studies have developed the Lanchester model for more realistic lethality rates \citep{taylor1975force,taylor1976canonical}; applied optimisation and game theoretic considerations \citep{taylor1977determining,taylor1979optimal,lin2014optimal,MacKay15}; and compared model outputs to the real world data sets of World War II \citep{lucas2004fitting,hung2005fitting}, the Vietnam War \citep{Schaffer68}, and the multi-adversary complexity of the Syrian civil war \citep{kress2018attrition}. More recently, \citep{Kalloniatis20Lanch} proposed a fully networked Lanchester model and established the emergence of behaviour reminiscent in manoeuvre warfare theory. Building on the theme of adding heterogeneity to Eq.(\ref{lanch}), recently \citep{Keane11} and \citep{Gonzalez11} built on the work of \citep{Protopopescu89} and added a spatial \textit{reaction-diffusion} component to Lanchester dynamics, which enabled concepts such as movement, terrain and perception to be reflected in the model.

The \textit{predator-prey} model of Lotka (\citeyear{Lotka1925}) and Volterra (\citeyear{Volterra1928}) is another adversarial model this work takes inspiration from. The simplest two-species instance of the model is given by, 
\begin{equation}
    \dot{N} = \beta^{(N)}_g N - \beta^{(N)}_d N P, \;\;
    \dot{P} = \beta^{(P)}_g N P - \beta^{(P)}_d P,
\label{LV}
\end{equation}
describing the interaction between a prey-species $N$, and predator-species $P$. Since its original inception, the highly idealised system in Eq.(\ref{LV}) has been generalised to include more realistic predation \citep{Holling66}, a general number of interacting species \citep{Hening18a,Hening18c}, and stochastic terms to simulate environmental effects \citep{Arato03}. 


\subsection{Exploring the role of influence in adversarial settings}
\label{influence}
Of particular interest to this work is the role that \textit{influence} plays in adversarial settings. This work is inspired by past OR studies which have applied ODEs to model insurgencies and cooperative dynamics. For example, \citep{MacKay15} combined the Lanchester, Deitchman and Richardson models in an attempt to demonstrate that state powers cannot deal with insurgencies solely by force; \citep{Kress09,Kress14} generalised the Deitchman model to explore the role that intelligence plays in counter-insurgency operations; \citep{Syms15} combined the Lanchester and logistic equations to model regime changes arising from insurgencies; and \citep{McLennan-Smith19, McLennan-Smith20} generalised Lotka-Volterra models to include symbiotic interactions to model adversarial dynamics in the vicinity of a homogeneous non-combatant humanitarian agencies. 


This current work is situated within the framework of multi-party Lanchester type models, such as \citep{MacKay09,MacKay12,kress2018attrition,Kalloniatis20Lanch}. As alluded,
this approach is particularly relevant to
counter-insurgency modelling, though not 
wholly bound by that context. We construct two relatively simple models which offer a means to explore the mutual influence between combatants and non-combatants in warfare. This is achieved by building upon the Lanchester-inspired ODE approaches cited in this work, and beyond, by adding inhomogeneity via logistic-Lotka-Volterra interactions to the non-combatant population. In particular, we exploit the ecology lens which has been richly discussed for counter-insurgency
for many years, though typically in qualitative
analogical terms \citep{Drapeau2008,Johnson2008,Zimmerman2016,Kilcullen2020}. 

In contrast to
\citep{McLennan-Smith19,McLennan-Smith20}
which also exploits ecological mechanisms
to model complex conflict environments, the present work presents the context where the civilian population may influence either force in their combatant role, rather than as humanitarian contribution. Naturally, the modern battlefield has both elements, as evidenced in concepts such as the three-block war \citep{Dorn09}.



\subsection{Outline of the paper}
In section \ref{SEC2} we present two instances of the model, one which conserves the total Green population over time --- labelled the \textit{supporter} model --- and one which lacks the conservation of the Green population --- labelled the \textit{contributor} model. To each model we offer example trajectories which inform our expectations for model behaviours. In section \ref{SEC3} we perform parameter sweeps which expose a number of interesting behaviours demonstrated by each of the models and discuss the likely reasons behind these outcomes. We also perform semi-analytic approximations of the model in a number of regimes in order to understand mechanisms behind these behaviours. Finally in section \ref{CONC} we discuss various implications stemming from outcomes and behaviours witnessed from the model and flag the potential for further work.

\section{The models}
\label{SEC2}
The two models we present in the following may be seen as a
segregation of a more complex unified model, which we separate in order to understand the dynamics peculiar to each. The overall framework is that using an ecology lens
to understand counter-insurgency, such as
Figure 2 of \citep{Drapeau2008}. Here, 
a spectrum of discrete states between `government'
forces, here seen as Blue, and `insurgent' combatants, or Red is decomposed: government infrastructure that provides resource
into Blue, government supporters, a broad political space
occupied by the population, then underground supporters
of the insurgency, and an insurgent auxiliary that can be
drawn into the conflict as Red combatants.
Thus we consider two segmentations of this spectrum: one between supporters and combatants, and the other
between contributors and combatants. These each require quite different types of non-linearities with implications we seek to understand.

\subsection{Model I: Green supporters}
\label{suppintro}
Consider three populations, \textit{Blue, Green} and \textit{Red}, labeled $B, G$ and $R$ respectively. $B$ and $R$ are adversaries, detrimentally affecting each other akin to the \textit{directed-fire} Lanchester-type interaction. The population $G$, is composed of three sub-populations, $G = g + \gamma + \Gamma$. Sub-populations $g$ and $\gamma$ are Red and Blue \textit{supporters}, respectively. Sub-population $\Gamma$ maintains neutrality and support neither side.  

The Lanchester interaction between Blue and Red is given by,
\begin{eqnarray}
\begin{split}
\dot{B}= - k^{\textrm{\tiny{L}}}_{\textrm{\tiny{R}}} f(g,\gamma) R \Theta (B),\;\; \dot{R} =- k^{\textrm{\tiny{L}}}_{\textrm{\tiny{B}}} f(\gamma,g) B \Theta(R),\\
\textrm{where} \;\; f(x,y) = \frac{x}{y+1},
\end{split}
\label{Lanchester1}
\end{eqnarray}
with the lethality coefficients $\{k^{\textrm{\tiny{L}}}_{\textrm{\tiny{B}}}, k^{\textrm{\tiny{L}}}_{\textrm{\tiny{R}}} \} \in \mathbb{R}_+$. Both expressions in Eq.~(\ref{Lanchester1}) are modulated by the \textit{dimensionless} term $f$, which factors the impact of both \textit{supporter} populations, $g$ and $\gamma$. The first argument in $f$ is the size of the \textit{supporting} population, and serves to increase the supported force's combat effectiveness against their adversary. For example, this may represent moral support which enables application of maximal combat power, or the utility of the supporting population's infrastructure. Correspondingly, the second argument in $f$ is the size of the \textit{detracting} population, decreasing the aforementioned force's effectiveness against the adversary. The term `1' in the denominator of $f$ represents a `standing population' that the detractors wish to maintain. If replaced by a number greater or less than than unity, the detractor population finds it easier or harder to influence the outcome of combat, respectively. Obviously, a choice close to zero would potentially lead to pathological behaviour in the model due to the appearance of unrealistically large values. Furthermore, the quantity $\Theta$ in both expressions in Eq.~(\ref{Lanchester1}) is a sufficiently smooth numerically stable Heaviside-step-like function which ensures that the populations of Blue and Red do not become negative. This function requires two parameters: the steepness of the transition and the offset from zero at which an entity is deemed extinct.
For the remainder of this work we use
\begin{equation}
    \Theta(x) = \frac{1}{2}\left\{ \tanh \left[\alpha \left( x-\frac{4}{\alpha} \right) \right] +1 \right\},
\end{equation}
where $\alpha = 10^6$ sets the scale for both the steepness
and the offset, chosen sufficiently large for a sharp transition, and the offset, up to this scale, is set at `4'. Intuitively, for sharp transitions the resultant dynamics will not be acutely sensitive to small variations in this offset.

The dynamics amongst the sub-populations of Green are determined by, 
\begin{eqnarray}
\begin{split}
\dot{g}= g \Gamma R \left( 1 - \frac{R}{k^{\textrm{\tiny{C}}}_{\textrm{\tiny{R}}}} \Theta(g) \right),\\
\dot{\gamma} =\gamma \Gamma B \left(1 - \frac{B}{k^{\textrm{\tiny{C}}}_{\textrm{\tiny{B}}}} \Theta(\gamma) \right),\\
\dot{\Gamma} =- g \Gamma R \left( 1 - \frac{R}{k^{\textrm{\tiny{C}}}_{\textrm{\tiny{R}}}}\Theta(g) \right) - \gamma \Gamma B \left(1 - \frac{B}{k^{\textrm{\tiny{C}}}_{\textrm{\tiny{B}}}} \Theta(\gamma) \right),
\end{split}
\label{Greenpop1}
\end{eqnarray}
where the inclusion of the $\Theta$ Heaviside terms ensure that the Green sub-populations remain positive. We remark that there is conservation of Green's total population \textit{i.e.} $\dot{G} \equiv \dot{g} + \dot{\gamma} + \dot{\Gamma} = 0$. 

The parameters $\{ k^{\textrm{\tiny{C}}}_{\textrm{\tiny{B}}}, k^{\textrm{\tiny{C}}}_{\textrm{\tiny{R}}}\} \in \mathbb{R}_+$ in Eq.~(\ref{Greenpop1}) are \textit{carrying capacities} associated with the corresponding logistic terms. The logistic terms are the mechanisms employed to increase/decrease the sub-population sizes. 
Here the carrying capacity plays
a significant role, well recognised
in insurgency modelling \citep{Drapeau2008}.
However, the limitations of one side or the
other to `carry' support go beyond
purely financial \citep{Syms15} or
terrain \citep{Johnson2008} factors, nor even
humanitarian \citep{Williamson2011} but also socio-cultural
\citep{Metz2004}. This includes the social
and tactical discipline of the respective
combatants. Thus, if the supported force ($R$ or $B$) becomes larger than the carrying capacity ($k^{\textrm{\tiny{C}}}_{\textrm{\tiny{R}}}$ or $k^{\textrm{\tiny{C}}}_{\textrm{\tiny{B}}}$) of their respective supporting Green sub-population ($g$ or $\gamma$), then the rate of change of the sub-population will become negative. This signals a rescinding of support, interpreted as non-combatants responding to negative consequences of military forces engaging in close proximity. Such consequences include \textit{collateral damage}, or \textit{cultural misunderstanding}. Thus a higher carrying capacity value ($k^{\textrm{\tiny{C}}}_{\textrm{\tiny{R}}}$ and $k^{\textrm{\tiny{C}}}_{\textrm{\tiny{B}}}$) corresponds to a more tolerant non-combatant population, and/or a military force which prioritises the best interests of the non-combatant population. 
The classic example of this
in 21st century experience is the Anbar Awakening
which presaged the 2007 surge
of US forces in Iraq \citep{Zimmerman2016}. These forces would possess adequate training, with doctrine which 
encourages military leaders to engage with community
leaders and facilitate their role in decision-making
\citep{Metz2004,Petraeus2009}. Evidently, different contexts
will have significantly different population carrying
capacities for a Blue force, contrasting the US occupations, respectively, of
Germany after World War II,
and Iraq after the Second Gulf War. 

In (Kress and Szechtman 2009) the authors proposed a Lanchester-insurgency model which considered the effect of a homogeneous and static Green population. The inclusion of the logistic terms in Eq.(6) enable a generalisation of the model in (Kress and Szechtman 2009) by considering the effects of a non-homogeneous Green which can both offer and rescind support to both forces.




\begin{figure}[!htb]
    \centering
    \includegraphics[width=12.9cm]{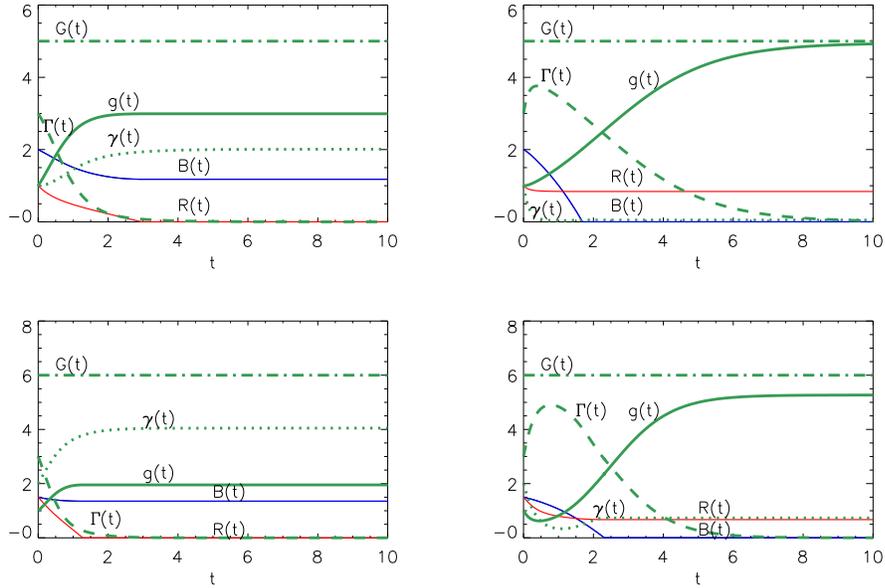}
    \caption{Examples of solutions for Model I. The red and blue curves correspond to $R$ and $B$ populations, respectively. The dotted green curve corresponds to $\gamma$ (support for $B$), the solid curve - to $g$ (support for $R$), the dashed green curve - to the neutral population $\Gamma$, and the dash-dotted green curve to the total $G$. The initial conditions are $B(0)=2$, $R(0)=1$, $g(0)=1$, $\gamma(0)=1$, $\Gamma(0)=3$ and $B(0)=1.5$, $R(0)=1.5$, $g(0)=1$, $\gamma(0)=2$, $\Gamma(0)=3$, for the top and bottom rows, respectively. The lethality coefficients are always $k^{\textrm{\tiny{L}}}_{\textrm{\tiny{R}}}=1$ and $k^{\textrm{\tiny{L}}}_{\textrm{\tiny{B}}}=1$. The carrying capacities are: $k^{\textrm{\tiny{C}}}_{\textrm{\tiny{R}}}=2$, $k^{\textrm{\tiny{C}}}_{\textrm{\tiny{B}}}=2$ (left column) and $k^{\textrm{\tiny{C}}}_{\textrm{\tiny{R}}}=1$, $k^{\textrm{\tiny{C}}}_{\textrm{\tiny{B}}}=1$ (right column).}
    \label{fig:fig_case1}
\end{figure}

Examples of solutions for Model I are shown in Figure~\ref{fig:fig_case1}. For all cases, demonstrated in the figure, despite the lethality coefficients being the same for Blue and Red populations, the outcome is not related to the initial populations, but strongly depends on the carrying capacity values. Namely, for the initially different Blue and Red populations (Blue is greater than Red; top row in Figure~\ref{fig:fig_case1}), Red wins for $k^{\textrm{\tiny{C}}}_{\textrm{\tiny{R}}}=1$ and $k^{\textrm{\tiny{C}}}_{\textrm{\tiny{B}}}=1$, while for larger capacities $k^{\textrm{\tiny{C}}}_{\textrm{\tiny{R}}}=2$ and $k^{\textrm{\tiny{C}}}_{\textrm{\tiny{B}}}=2$, Blue wins. This is due to the support offered by Green (solid green curve). 
\subsection{Model II: Green contributors}
\label{contribintro}
In the second model instance, the sub-populations $g$ and $\gamma$ play the role of combat \textit{contributors} to Red and Blue, respectively. As with the previous \textit{supporter} model, sub-population $\Gamma$ maintains neutrality. The Lanchester-type interaction between Blue and Red is given by,
\begin{eqnarray}
\begin{split}
    \dot{B} = - R \Theta(B) + k^{\textrm{\tiny{T}}}_{\textrm{\tiny{B}}} \gamma \Theta(\gamma), \;\; \dot{R} = - B \Theta(R) + k^{\textrm{\tiny{T}}}_{\textrm{\tiny{R}}} g \Theta(g).
\end{split}
\label{Lanchester2}
\end{eqnarray}
Eq. (\ref{Lanchester2}) differs from the interactions presented in Eq.~(\ref{Lanchester1}) by placing the sub-populations of $\gamma$ and $g$ on the \textit{same level} as the adversaries $B$ and $R$, rather than as arguments into the modulating the dimensionless factor $f$. Sub-populations $\gamma$ and $g$ now act as sources, and are employed to solely decrease the effectiveness of forces $B$ and $R$, respectively. Thus, the meaning of the parameters $\{k^{\textrm{\tiny{T}}}_{\textrm{\tiny{B}}},k^{\textrm{\tiny{T}}}_{\textrm{\tiny{R}}}\}$ has changed from the previous case, being referred to Blue and Red's \textit{transfer} coefficients, respectively. Additionally, compared to Eq.(\ref{Lanchester1}), we have set the lethality coefficients to unity in this instance, which can be generalised in the future. Again, the inclusion of the multiplicative $\Theta$ terms ensures the positivity of population values. In order to reflect the different way that $\gamma$ and $g$ interact in the adversarial setting in Eq.~(\ref{Lanchester2}), we change the dynamics amongst the three sub-populations of Green to be non-conservative via,
\begin{eqnarray}
\begin{split}
   \dot{g} = g \Gamma R \left( 1 - \frac{R}{k^{\textrm{\tiny{C}}}_{\textrm{\tiny{R}}}} \Theta(g)\right) - k^{\textrm{\tiny{T}}}_{\textrm{\tiny{R}}} g \Theta(g) ,\\
   \dot{\gamma} = \gamma \Gamma B \left(1 -\frac{B}{k^{\textrm{\tiny{C}}}_{\textrm{\tiny{B}}}}\Theta(\gamma) \right) - k^{\textrm{\tiny{T}}}_{\textrm{\tiny{B}}} \gamma \Theta(\gamma),\\
   \dot{\Gamma} = - g \Gamma R \left( 1 - \frac{R}{k^{\textrm{\tiny{C}}}_{\textrm{\tiny{R}}}}\Theta(g) \right) - \gamma \Gamma B \left( 1 - \frac{B}{k^{\textrm{\tiny{C}}}_{\textrm{\tiny{B}}}}\Theta(\gamma) \right).
\end{split}
\label{Greenpop2}
\end{eqnarray}
As with Eq.~(\ref{Greenpop1}) the constants $\{ k^{\textrm{\tiny{C}}}_{\textrm{\tiny{B}}}, k^{\textrm{\tiny{C}}}_{\textrm{\tiny{R}}}\} \in \mathbb{R}_+ $ play the role of carrying capacities, with their function and implied meaning unchanged from the first model discussed in section \ref{suppintro}. Furthermore, the defining equations for $g$ and $\gamma$ now contain additional sink terms, balancing the source terms in the defining equations for Blue and Red. Thus, the total rate of change of $G$ is now given by $\dot{G} = - k^{\textrm{\tiny{T}}}_{\textrm{\tiny{B}}} \gamma \Theta(\gamma) - k^{\textrm{\tiny{T}}}_{\textrm{\tiny{R}}} g \Theta(g)$. 

\begin{figure}[!htb]
    \centering
    \includegraphics[width=12.9cm]{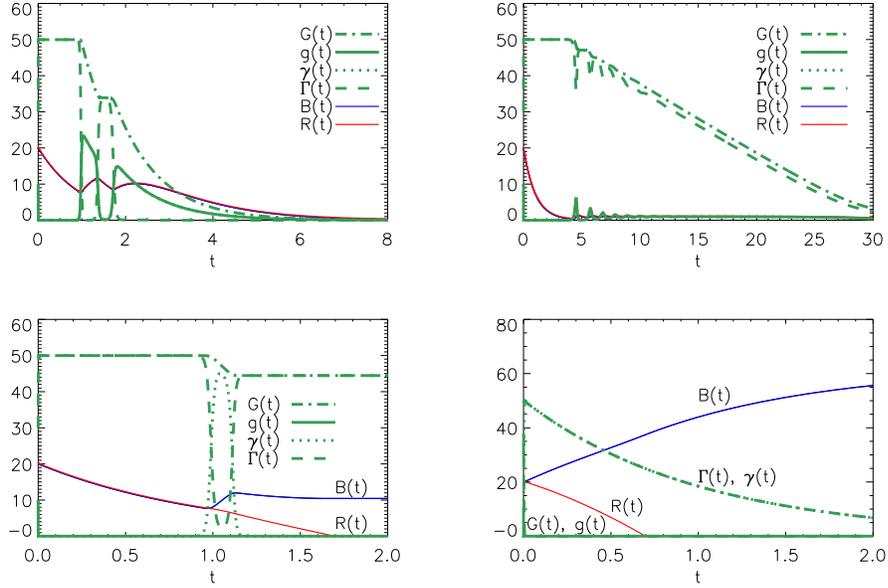}
    \caption{Example outputs of Model II for various parameter choices. Colour meanings for trajectories are same as in figure \ref{fig:fig_case1}. The initial conditions are $B(0)=20$, $R(0)=20$, $g(0)=10$, $\gamma(0)=10$, $\Gamma(0)=10$ and same for all given examples. The transfer coefficients are always $k^{\textrm{\tiny{T}}}_{\textrm{\tiny{B}}}=1$ and $k^{\textrm{\tiny{T}}}_{\textrm{\tiny{R}}}=1$. The carrying capacities are: $k^{\textrm{\tiny{C}}}_{\textrm{\tiny{R}}}=10$, $k^{\textrm{\tiny{C}}}_{\textrm{\tiny{B}}}=10$ (top-left plot); $k^{\textrm{\tiny{C}}}_{\textrm{\tiny{R}}}=1$, $k^{\textrm{\tiny{C}}}_{\textrm{\tiny{B}}}=1$ (top-right plot); $k^{\textrm{\tiny{C}}}_{\textrm{\tiny{R}}}=1$, $k^{\textrm{\tiny{C}}}_{\textrm{\tiny{B}}}=10$ (bottom-left plot); $k^{\textrm{\tiny{C}}}_{\textrm{\tiny{R}}}=10$, $k^{\textrm{\tiny{C}}}_{\textrm{\tiny{B}}}=100$ (bottom-right plot).}
    \label{fig:fig_case2}
\end{figure} 

This model instance focuses on the same initial conditions for populations $B$ and $R$, and $\gamma$ and $g$, in order to understand the effect that Green has on the outcome for near-peer adversaries. Additionally, varying initial conditions in the contributor model matters less as $B$ and $R$ can now grow in population, resulting in both populations potentially sampling more of the phase space for each trajectory. Future work concerning this model may address different initial conditions more thoroughly. Examples of solutions for this model instance are shown in Figure~\ref{fig:fig_case2}. Larger initial conditions were chosen to more easily demonstrate details in trajectories. For all panels the transfer coefficients are the same with $k^{\textrm{\tiny{T}}}_{\textrm{\tiny{B}}}= k^{\textrm{\tiny{T}}}_{\textrm{\tiny{R}}} = 1$. Outcomes and trajectories differ greatly however, strongly depending on the values of the carrying capacities. Focusing on the top row for equal carrying capacities the result is a draw; though this stalemate is achieved through an oscillatory stage, whose temporal extent depends on the carrying capacities. For large $k^{\textrm{\tiny{C}}}_{\textrm{\tiny{R}}}=k^{\textrm{\tiny{C}}}_{\textrm{\tiny{B}}}=10$ (top-left), the oscillatory stage is short, and the Green population decreases faster than the Blue or Red populations. For small $k^{\textrm{\tiny{C}}}_{\textrm{\tiny{R}}}=k^{\textrm{\tiny{C}}}_{\textrm{\tiny{B}}}=1$ (top-right), the Red and Blue populations quickly decrease to almost zero, while the Green population keeps its oscillatory decrease until exhausted. We will address the nature of oscillations in the solution in Section~\ref{SEC3}. 

For different carrying capacities, the outcome depends on whether one of the carrying capacities is greater than the corresponding Red or Blue population. For $k^{\textrm{\tiny{C}}}_{\textrm{\tiny{R}}}=1$ and $k^{\textrm{\tiny{C}}}_{\textrm{\tiny{B}}}=10$ (bottom-left; both carrying capacities are less than $B(0)$ and $R(0)$), the Red and Blue populations initially decrease until Blue reaches the threshold level. Then, Blue population receives support from the Green population, which consequently leads to Blue's victory. However, for $k^{\textrm{\tiny{C}}}_{\textrm{\tiny{R}}}=10$ and $k^{\textrm{\tiny{C}}}_{\textrm{\tiny{B}}}=100$ (bottom-right; $k^{\textrm{\tiny{C}}}_{\textrm{\tiny{R}}} > B(0)$), Blue gets immediate support from the Green population, which leads to Blue victory. The behaviour of the Green population is very different in the two latter cases. For small carrying capacities, the Green population remains constant before and after the critical phase at $t=1$, and a small decrease in Green population leads to the Blue victory. For large carrying capacities, Green contributes to Blue over the extent of the solution and keeps contributing to Blue even after Red population reaches zero at $t=0.7$.

\section{Numerical analysis and semi-analytic approximations}
\label{SEC3}
In order to understand model behaviours for both model iterations, we offer more substantial numerical analysis which highlight differences due to various parameters changes. Furthermore, in key instances we \textit{solve} appropriate model approximations which enable a more in-depth understanding of interesting model behaviours. 


\subsection{Supporting case: parameter sweeps}
\label{suppsweep}
\begin{figure}[!htb]
    \centering
    \includegraphics[width=13.9cm]{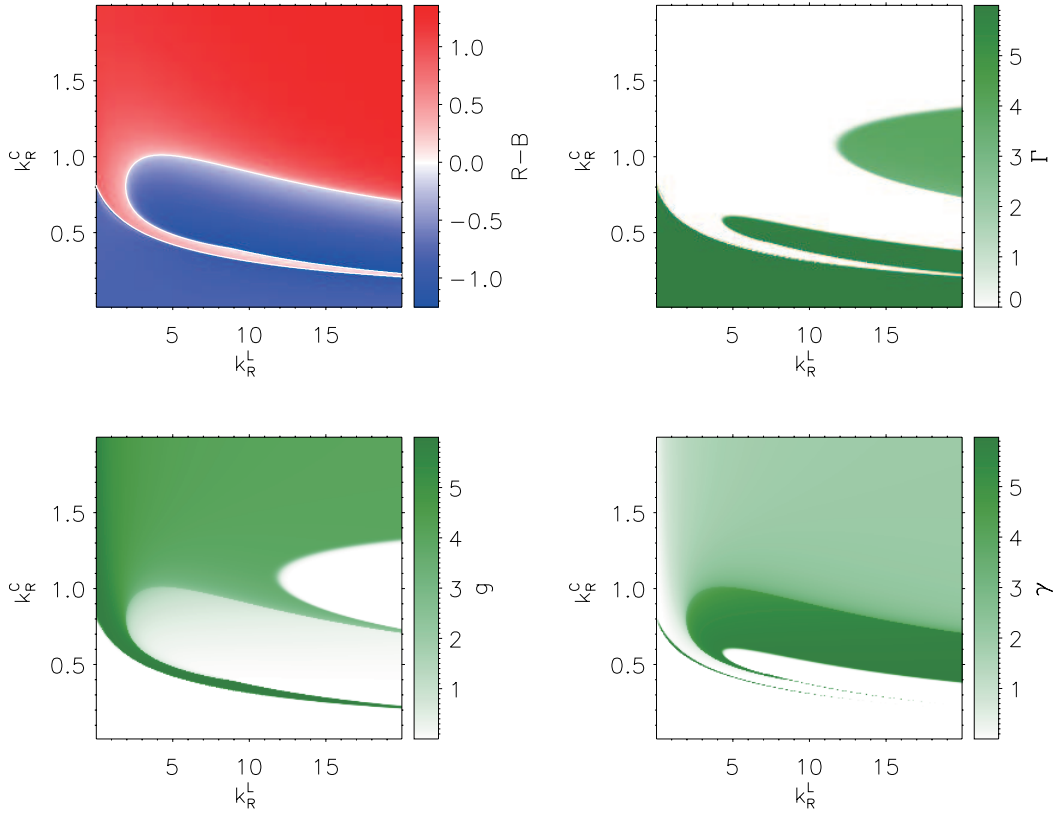}
    \caption{Parameter scans for Model I, $B(t=0)=1.5$, $R(t=0)=1.5$, $g(t=0)=1$, $\gamma(t=0)=2$, $\Gamma(t=0)=3$, $k^{\textrm{\tiny{L}}}_{\textrm{\tiny{B}}} = 1$, $k^{\textrm{\tiny{C}}}_{\textrm{\tiny{B}}} = 1$. Colours indicate the final population values (at $t=50$) over a range of $k^{\textrm{\tiny{L}}}_{\textrm{\tiny{R}}}$ and $k^{\textrm{\tiny{C}}}_{\textrm{\tiny{R}}}$, with blue and red signifying $B$ and $R$ victories in the top left panel. The white contour in the top-left panel separates the R and B victory regions.}
    \label{fig:fig_case1_scan}
\end{figure}

For the Green-supporting scenario given by Eqs.(\ref{Lanchester1},\ref{Greenpop1}), we offer phase plots in Figure \ref{fig:fig_case1_scan} which detail final values (taken at $t=50$) for the populations over a range of $k^{\textrm{\tiny{L}}}_{\textrm{\tiny{R}}} \in (0,20)$ (Red's lethality coefficient) and $k^{\textrm{\tiny{C}}}_{\textrm{\tiny{R}}} \in (0,2)$ (carrying capacity dictating $g$'s capacity to support Red). Thus, in Figure \ref{fig:fig_case1_scan} we are examining model outcomes of Red changing its parameter inputs --- combat effectiveness, and it's ability to maintain support from $g$. In all four panels the initial conditions are $B_0 = R_0 = 1.5$, $g_0 =1$, $\gamma_0 =2$ and $\Gamma_0 = 3$.

The top left panel of Figure \ref{fig:fig_case1_scan} details the final value of the expression $R-B$ for Eq.(\ref{Lanchester1}). Positive values (Red wins) are indicated by shades of red, and negative values (Blue wins) indicated by shades of blue. White contours separate the regions between Red and Blue victories. In the region $k^{\textrm{\tiny{C}}}_{\textrm{\tiny{R}}} > 1$ Red always wins engagements. Correspondingly, for low values of $k^{\textrm{\tiny{C}}}_{\textrm{\tiny{R}}} < 0.2$ Blue is always the victor. In the region $0.2< k^{\textrm{\tiny{C}}}_{\textrm{\tiny{R}}} <1$ the outcome is more complicated due to the presence of two white contours (signifying a stalemate) separating the regions of Red and Blue winning. Thus in this region there is a thin slice between the majority of Blue victories where Red emerges victorious. 

Contrast to this, on the bottom-left panel for Red contributors $g$, the regions where $g$ are non-zero correspond to the regions where Red wins the engagement in the top left diagram. This includes the appearance of the thin slice between the Blue victories in the region $k^{\textrm{\tiny{C}}}_{\textrm{\tiny{R}}} \in (0.2, 1)$. This slice also exists as non-zero $g$ values in the bottom-left panel. Nevertheless, there is a notable half-oval shape in the parameter values $k^{\textrm{\tiny{L}}}_{\textrm{\tiny{R}}} \times k^{\textrm{\tiny{C}}}_{\textrm{\tiny{R}}} \in (10,20) \times ( 1,1.5)$, absent from the top-left diagram, indicating a region of zero $g$-values, even though Red wins the engagement in this region. Trajectories in this region reveal that Blue reaches zero population before Red descends below $k^{\textrm{\tiny{C}}}_{\textrm{\tiny{R}}}$-value --- \textit{i.e.} Red \textit{does not} have the ability to sustain support from Green. Thus in this region $g$'s population declines after the engagement ends. Just below this region in the $k^{\textrm{\tiny{C}}}_{\textrm{\tiny{R}}}$ direction Red still obtains victory, but the final values of $g$ are non-zero. Trajectories in this region show that although $g$-values initially begin decreasing, once Blue loses the engagement Red population is below $k^{\textrm{\tiny{C}}}_{\textrm{\tiny{R}}}$-value, meaning that $g$ increases after the engagement has ended. Decreasing $k^{\textrm{\tiny{C}}}_{\textrm{\tiny{R}}}$ values there is an enclosed region where Blue obtains victory, in addition to small (or zero) $g$-values at the final time. Trajectories in this region reveal that although Blue starts the engagement against Red poorly, because $k^{\textrm{\tiny{C}}}_{\textrm{\tiny{R}}}<k^{\textrm{\tiny{C}}}_{\textrm{\tiny{B}}}$, and hence Red's ability to sustain Green support is less than Blue's, Blue eventually wins the engagement due to experiencing much more support. 

Focusing on lower $k^{\textrm{\tiny{C}}}_{\textrm{\tiny{R}}}$ values, we reach the aforementioned strip of red surrounded by Blue victories. In this region Red wins the engagement in combination with large $g$-values at the final time. Trajectories in this region reveal that Blue is actually ahead in engagements for the majority of the time. However, as Red's population value for the majority of time is less than $k^{\textrm{\tiny{C}}}_{\textrm{\tiny{R}}}$, whereas Blue's is greater than $k^{\textrm{\tiny{C}}}_{\textrm{\tiny{B}}}$, Red is able to maintain support from $g$ for a significant amount of time, whereas Blue's support from $\gamma$ dwindles away. Red's growing support from Green means that Blue's combat effectiveness eventually plummets, leading to Red victory. In the region for lowest $k^{\textrm{\tiny{C}}}_{\textrm{\tiny{R}}}$ values which show Blue victory in the top panel, there are zero $g$ \textit{and} $\gamma$ values. In this region, both Red and Blue \textit{completely} lose support from their respective Green populations. Nevertheless, due to Blue's ability to sustain Green support being greater than Red's, Red has less (though non-zero) force than Blue by the final time.

\subsection{Supporting case: analytic considerations}
\label{modIanalyt}
Using approximations it is possible to gain an analytical understanding of the mechanisms causing the bottom white contour in the top-left plot of Figure \ref{fig:fig_case1_scan}, below which Green has completely rescinded support from both Red and Blue. To enable this we assign $k^{\textrm{\tiny{L}}}_{\textrm{\tiny{R}}}=k^{\textrm{\tiny{L}}}_{\textrm{\tiny{B}}}$ and $k^{\textrm{\tiny{C}}}_{\textrm{\tiny{R}}} = k^{\textrm{\tiny{C}}}_{\textrm{\tiny{B}}}$ in Eq.(\ref{Greenpop1}), and consider the case $\gamma = g$, thus focusing on the case where Green supports both forces equally. Furthermore, if we assume that the starting populations of Blue and Red are equal (\textit{i.e} $B_0 = R_0$), then this forces $B=R$ in Eq.(\ref{Lanchester1}). 
The significance of this assumption is that
we focus on what may be deemed the worst-case scenario
of stalemate where conflict is drawn out with
devastating and exhausting impacts on the people
and resources of both sides of the conflict.

Thus the system under consideration becomes
\begin{equation}
    \dot{B} = - k^{\textrm{\tiny{L}}}_{\textrm{\tiny{R}}} \frac{g}{1+g} B, \;\;
    \dot{g} = g (G_0 - 2g ) B \left( 1 - \frac{B}{k^{\textrm{\tiny{C}}}_{\textrm{\tiny{R}}}} \right),
\label{sys1}
\end{equation}
where the initial condition for $g$ is bounded by $0 < g_0 < G_0/2$. Eq.(\ref{sys1}) reveals two important facts:
\begin{itemize}
    \item{the right hand side of the equation for $\dot{B}$ is always less than or equal to zero}
    \item{if $B_0 \le k^{\textrm{\tiny{C}}}_{\textrm{\tiny{R}}}$, then the right hand side of the equation for $g$ will remain positive, ensuring that $B \rightarrow 0$ after sufficient time.}
\end{itemize} 


We solve Eq.(\ref{sys1}) to obtain the relationship between dynamic variables $B$ and $g$ and determine the remaining conditions for $B_0 > k^{\textrm{\tiny{C}}}_{\textrm{\tiny{R}}}$ which enable $B \rightarrow 0$ after sufficient time. Dividing $\dot{B}$ by $\dot{g}$ in Eq.(\ref{sys1}) reveals
\begin{eqnarray}
\begin{split}
    &\frac{d B}{d g} = -\frac{k^{\textrm{\tiny{L}}}_{\textrm{\tiny{R}}}}{(1+g)(G_0-2 g) \left(1 - \frac{B}{k^{\textrm{\tiny{C}}}_{\textrm{\tiny{R}}}}\right)} ,\\
    \Rightarrow \;\;& (B-B_0)(B+B_0 -2 k^{\textrm{\tiny{C}}}_{\textrm{\tiny{R}}}) \\
    &= \frac{2 k^{\textrm{\tiny{L}}}_{\textrm{\tiny{R}}} k^{\textrm{\tiny{C}}}_{\textrm{\tiny{R}}}}{G_0+2} \ln \left[\left(\frac{1+g}{G_0 - 2g}\right)\left(\frac{G_0-2g_0}{1+g_0}\right)\right].
\end{split}
\label{Bofg}
\end{eqnarray}
Expressing $g$ as a function of $B$, Eq.(\ref{Bofg}) becomes
\begin{eqnarray}
\begin{split}
   & g = \frac{G_0 h(B) -1}{1+2h(B)}, \\
   \textrm{where} \;\;& h(B) = \frac{1+g_0}{G_0-2g_0} e^{  \frac{G_0+2}{2 k^{\textrm{\tiny{L}}}_{\textrm{\tiny{R}}} k^{\textrm{\tiny{C}}}_{\textrm{\tiny{R}}}}(B-B_0)(B+B_0 -2 k^{\textrm{\tiny{C}}}_{\textrm{\tiny{R}}})}.
\end{split}
    \label{gofB}
\end{eqnarray}
The form of the equation for $B$ in Eq.(\ref{sys1}) means that the magnitude of $B$ can only ever decrease from the initial value of $B_0$. However, using Eq.(\ref{gofB}), we can determine the values of the initial conditions, and $\{k^{\textrm{\tiny{L}}}_{\textrm{\tiny{R}}}, k^{\textrm{\tiny{C}}}_{\textrm{\tiny{R}}}\}$ which allow $B$ to reach a value of zero. Due to the numerator of $g(B)$ in Eq.(\ref{gofB}), the Green sub-populations have the potential to be completely depleted if the minimum value of $h(B) \le 1/ G_0$. As the argument of $h(B)$ in Eq.(\ref{gofB}) is a positive quadratic, with minimum value $B=k^{\textrm{\tiny{C}}}_{\textrm{\tiny{R}}}$, the minimum value of $h(B)$ is $h(k^{\textrm{\tiny{C}}}_{\textrm{\tiny{R}}})$. Thus, the condition between the initial conditions and $\{ k^{\textrm{\tiny{L}}}_{\textrm{\tiny{R}}}, k^{\textrm{\tiny{C}}}_{\textrm{\tiny{R}}}\}$ which enables $B \rightarrow 0$ is given by
\begin{eqnarray}
\begin{split}
& h(k^{\textrm{\tiny{C}}}_{\textrm{\tiny{R}}}) > \frac{1}{G_0},\\
    \Rightarrow \;\;& k^{\textrm{\tiny{C}}}_{\textrm{\tiny{R}}} > B_0 +  \rho(k^{\textrm{\tiny{L}}}_{\textrm{\tiny{R}}}) - \sqrt{ \rho(k^{\textrm{\tiny{L}}}_{\textrm{\tiny{R}}}) [ \rho(k^{\textrm{\tiny{L}}}_{\textrm{\tiny{R}}}) + 2 B_0]},
    \end{split}
\label{k1andk3}
\end{eqnarray}
where
\begin{equation}
    \rho(k^{\textrm{\tiny{L}}}_{\textrm{\tiny{R}}}) = \frac{k^{\textrm{\tiny{L}}}_{\textrm{\tiny{R}}}}{G_0+2} \ln \frac{G_0 (1+g_0)}{G_0-2 g_0}.
\end{equation}
Thus, given initial conditions $G_0$, $g_0 < G_0/2$ and $B_0 > k^{\textrm{\tiny{C}}}_{\textrm{\tiny{R}}}$, the values of $k^{\textrm{\tiny{L}}}_{\textrm{\tiny{R}}}$ and $k^{\textrm{\tiny{C}}}_{\textrm{\tiny{R}}}$ which satisfy the requirement given by Eq.(\ref{k1andk3}) will ensure that $B \rightarrow 0$ after sufficient time.

\begin{figure}
\begin{center}
\includegraphics[width=7cm]{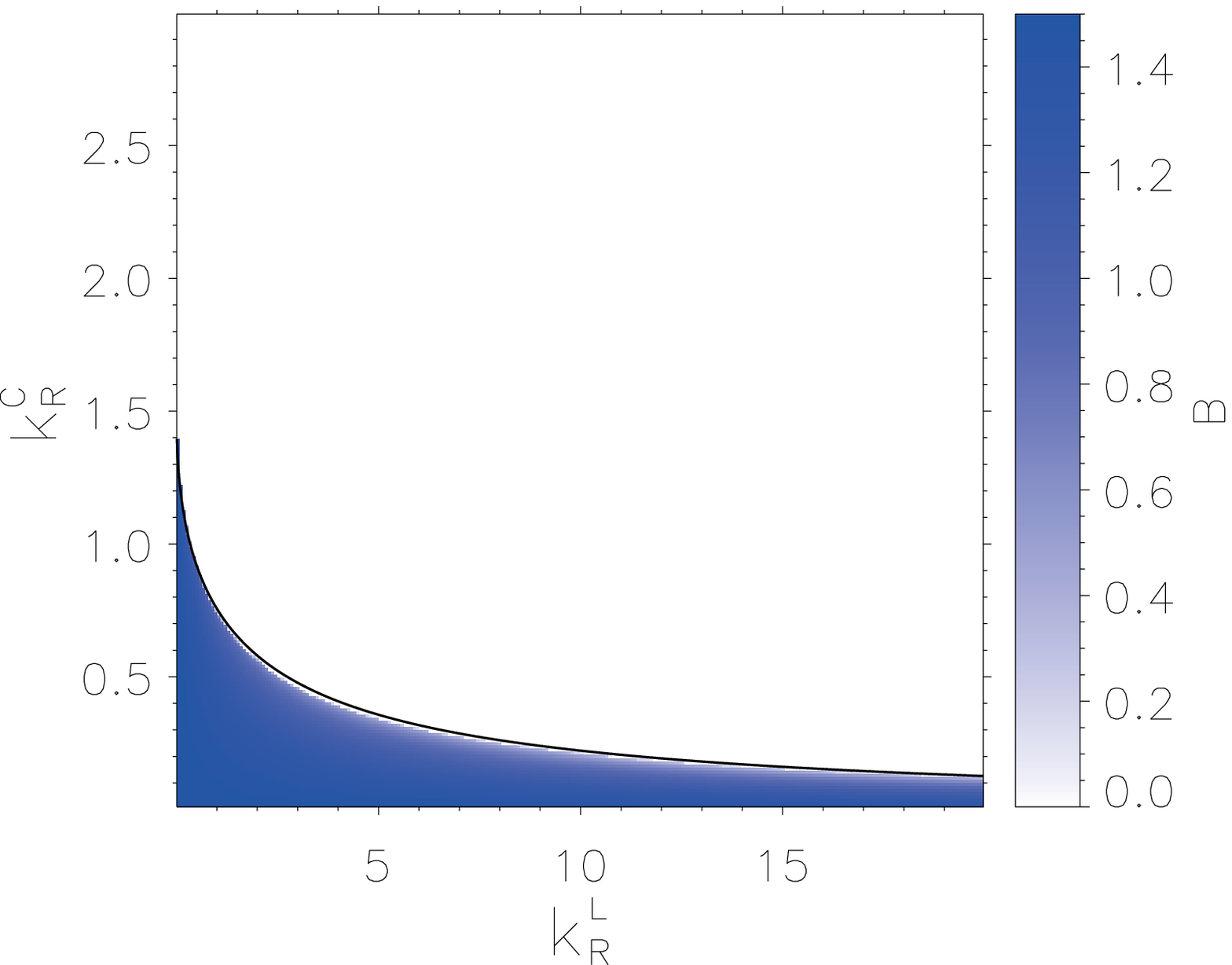}
\caption{Comparison of Eq.(\ref{k1andk3}) with corresponding numerical outputs of final $B$-values of Eqs.(\ref{Lanchester1},\ref{Greenpop1}) with $k^{\textrm{\tiny{L}}}_{\textrm{\tiny{R}}} = k^{\textrm{\tiny{L}}}_{\textrm{\tiny{B}}}$, $k^{\textrm{\tiny{C}}}_{\textrm{\tiny{R}}} = k^{\textrm{\tiny{C}}}_{\textrm{\tiny{B}}}$. Initial conditions applied are $B_0=R_0=1.5$, $g_0=\gamma_0=1$, $\Gamma_0=3$, resulting in system becoming Eq.(\ref{sys1}) --- \textit{i.e.} $B=R$ and $g=\gamma$.  The black contour gives the case of Eq.(\ref{k1andk3}) with an equals sign which is the boundary between zero and non-zero $B$ values at the final time.}
\label{fig:comparwitk1andk3}
\end{center}
\end{figure}

Figure \ref{fig:comparwitk1andk3} reveals very good agreement of Eq.(\ref{k1andk3}) with corresponding numerical outputs of final $B$-values stemming from  Eqs.(\ref{Lanchester1},\ref{Greenpop1}) with $k^{\textrm{\tiny{L}}}_{\textrm{\tiny{R}}} = k^{\textrm{\tiny{L}}}_{\textrm{\tiny{B}}}$, $k^{\textrm{\tiny{C}}}_{\textrm{\tiny{R}}} = k^{\textrm{\tiny{C}}}_{\textrm{\tiny{B}}}$, and matching initial conditions for $B,R$ and $g, \gamma$. Eq.(\ref{k1andk3}) allows for effective interrogation of the impact of changes to parameters and initial conditions on the form of the contour which delineates between Green contributing, or choosing to withdraw. For instance, if we mandate that $|2B_0/\rho(k^{\textrm{\tiny{L}}}_{\textrm{\tiny{R}}})|<1$, then performing a binomial expansion on the square root of Eq.(\ref{k1andk3}) obtains
\begin{equation}
    k^{\textrm{\tiny{C}}}_{\textrm{\tiny{R}}} > \frac{B^2_0 (G_0+2)}{2 k^{\textrm{\tiny{L}}}_{\textrm{\tiny{R}}} \ln \frac{G_0 (1+g_0)}{G_0-2 g_0}} + O \left( \frac{B^3_0}{(k^{\textrm{\tiny{L}}}_{\textrm{\tiny{R}}})^2} \right).
    \label{carcapreq}
\end{equation}
Eq.(\ref{carcapreq}) reveals the relationships between carrying capacities and lethalities for peer forces pursuing mutual annihilation with support from from the non-combatant population --- a situation to be avoided.   

\subsection{Contributing case: parameter sweeps}
\label{contribsweep}
\begin{figure*}[!htb]
    \centering
    \includegraphics[width=14.3cm]{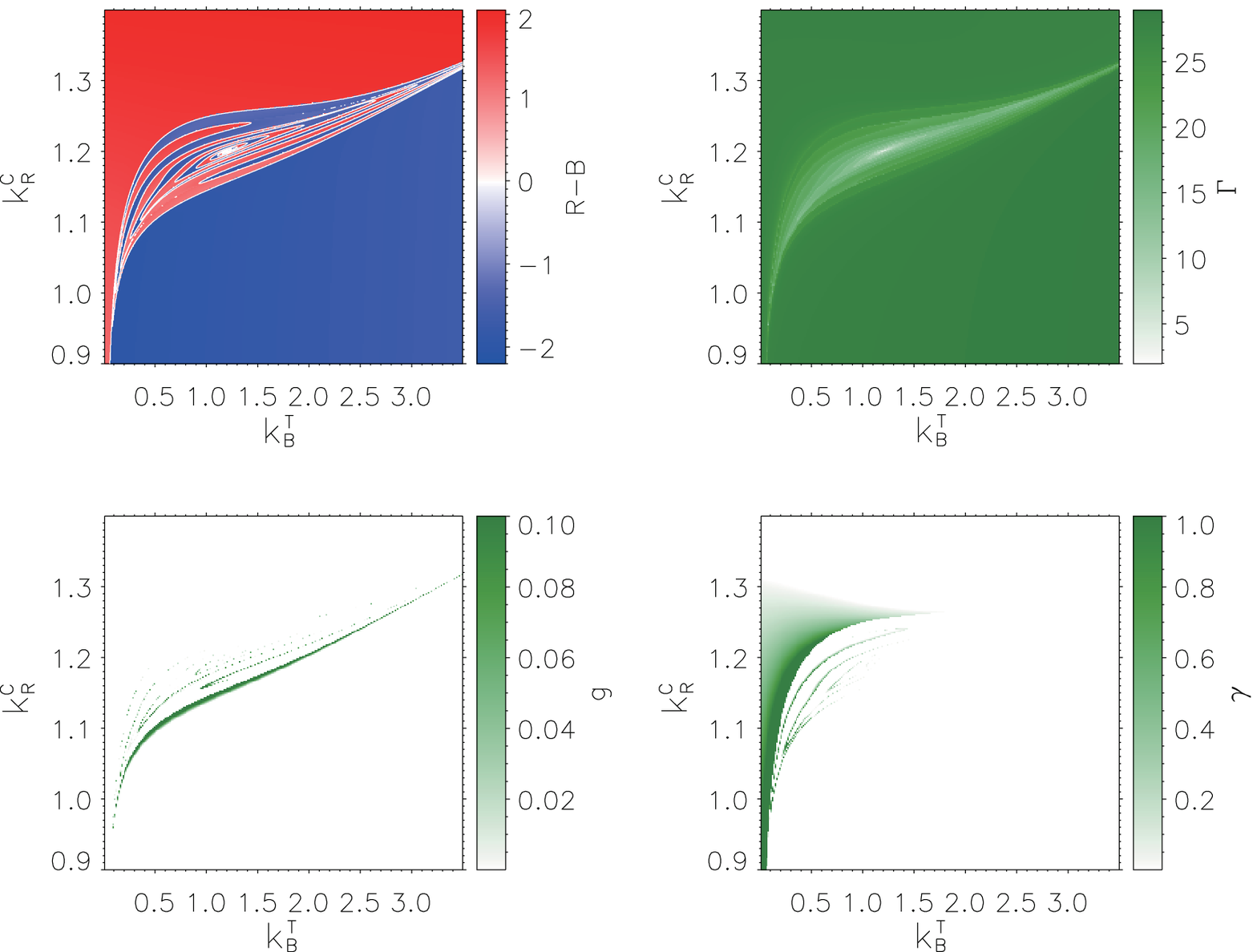}
    \caption{Parameter scans for Model II. The solutions are produced for $B(t=0)=R(t=0)=20$, $g(t=0)=\gamma(t=0)=\Gamma(t=0)=10$, and $k^{\textrm{\tiny{T}}}_{\textrm{\tiny{R}}}=k^{\textrm{\tiny{C}}}_{\textrm{\tiny{B}}}=1.2$. Colours indicate the final population values (at $t=50$) over a range of $k^{\textrm{\tiny{T}}}_{\textrm{\tiny{B}}}$ and $k^{\textrm{\tiny{C}}}_{\textrm{\tiny{R}}}$, with blue and red signifying $B$ and $R$ victories in the top left panel. The white regions in the top-left panel separate the Red and Blue winning regions.}
    \label{fig:fig_case2_scan}
\end{figure*}

We offer phase plots of final population values for the Green-contributing scenario in Figure \ref{fig:fig_case2_scan}, over the range $k^{\textrm{\tiny{T}}}_{\textrm{\tiny{B}}} \in (0,3.5)$ (Blue's transfer coefficient) and $k^{\textrm{\tiny{C}}}_{\textrm{\tiny{R}}} \in (0.9,1.4)$ (Red's capacity to maintain support from $g$). A narrower parameter range than those considered in Figure \ref{fig:fig_case1_scan} is given due to the contours being significantly more complicated in the contributor scenario. The top-left panel for $R-B$ values at $t=50$ generally has three distinct regions: top-left diagonal of Red victories; bottom-right diagonal of Blue victories; and middle-diagonal consisting of a complex region of mixed outcomes separated by a highly deformed spiral. 

In the top-left and bottom-right diagonal regions, trajectories generally follow the pattern of Green contributors initially retreating to the neutral $\Gamma$ population due to Red and Blue being over their respective carrying capacity values. Thus Red and Blue initially damage each other at the same rate as their starting populations are equal. At some point in time in either of these regions, either the Red or Blue population becomes lower than their carrying capacity value. Once this happens, the following outcomes are possible:
\begin{itemize}
    \item{In the upper-left diagonal region with $k^{\textrm{\tiny{C}}}_{\textrm{\tiny{R}}} > k^{\textrm{\tiny{C}}}_{\textrm{\tiny{B}}} $, above the white contours, Red's capacity to maintain support from Green is sufficiently greater than Blue's, and Blue's transfer coefficient ($k^{\textrm{\tiny{T}}}_{\textrm{\tiny{B}}}$) is sufficiently less than Red's. This leads to a sharp rise in Green population $g$ who contribute to Red's effort, which in turn leads to a rise in Red's population which is enough to secure victory against Blue before any support for Blue can take effect.}
    \item{For $k^{\textrm{\tiny{C}}}_{\textrm{\tiny{R}}} < k^{\textrm{\tiny{C}}}_{\textrm{\tiny{B}}} $ in the same upper-left diagonal region, Blue's population decreases more rapidly, leading to $\gamma$ increasing in population relatively early in the engagement. However, because Blue's transfer coefficient $k^{\textrm{\tiny{T}}}_{\textrm{\tiny{B}}}$ is so small, the majority of the corresponding contributors from Green never get the chance to \textit{actually} contribute. Hence, Red emerge victorious as they are much more effective at transferring their Green contributors into the engagement after they start receiving Green support.}
    \item{In the lower-right diagonal region with $k^{\textrm{\tiny{C}}}_{\textrm{\tiny{R}}} > k^{\textrm{\tiny{C}}}_{\textrm{\tiny{B}}} $, below the white contours, Red's capacity to obtain contribution from Green is greater than Blue's. Thus Red's contributor population from Green begins increasing before Blue's. Nevertheless in this region Blue's transfer coefficient ($k^{\textrm{\tiny{T}}}_{\textrm{\tiny{B}}})$ is sufficiently large enough that when Blue finally does obtain contributions from Green, Blue's ability to transfer said contributors into the engagement is enough to counter Red's efforts and win the engagement --- even though Blue's action was slower compared to Red's.}
    \item{For $k^{\textrm{\tiny{C}}}_{\textrm{\tiny{R}}} < k^{\textrm{\tiny{C}}}_{\textrm{\tiny{B}}}$ in the lower-right diagonal region, as with the equivalent upper-left diagonal region, Blue reaches the capacity to obtain contributions from Green before Red, meaning that $\gamma$ increases in population before $g$. Blue's transfer coefficient $k^{\textrm{\tiny{T}}}_{\textrm{\tiny{B}}}$ is also sufficiently large that the Green contributors are transferred effectively to swiftly ensure Blue's victory.}
\end{itemize}
Generally in these regions above and below the white contours, the initial and final total Green populations do not differ substantially from each other --- as seen by the top-right panel in Figure \ref{fig:fig_case2_scan} for final $\Gamma$ population which does not differ greatly from $G_0$ in these regions.

\begin{figure*}[!htb]
\begin{center}
\includegraphics[width=14.3cm]{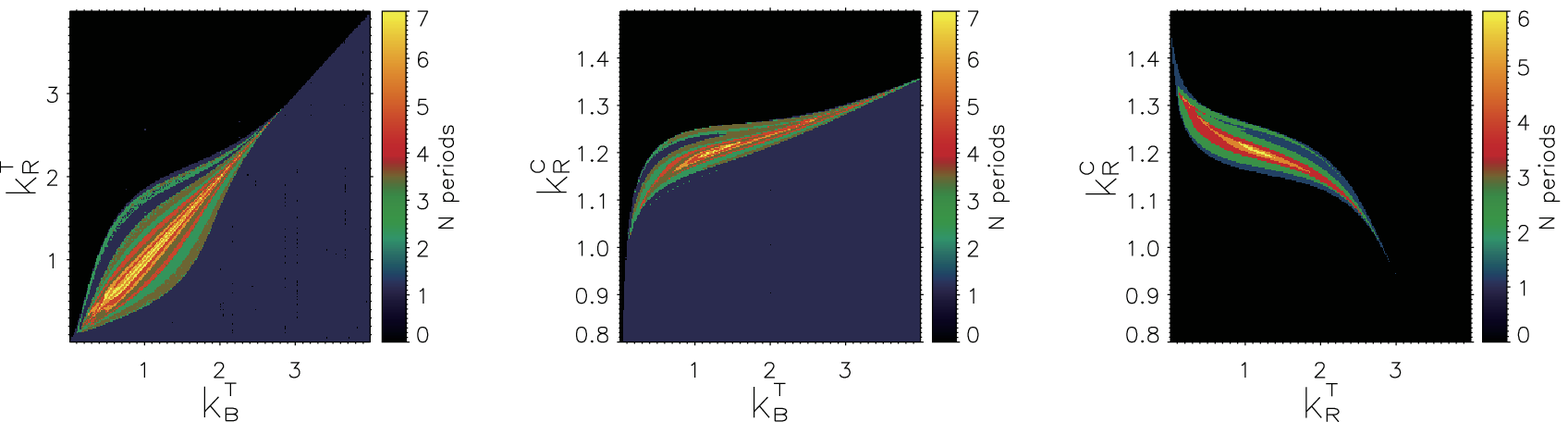}
\caption{Number of periods in $R$ --- signified by different colours ranging from $0$ (black) to $7$ (yellow) periods. Numerical outputs are produced for $B(t=0)=R(t=0)=20$, $g(t=0)=\gamma(t=0)=\Gamma(t=0)=10$. Left panel: $k^{\textrm{\tiny{C}}}_{\textrm{\tiny{R}}}=k^{\textrm{\tiny{C}}}_{\textrm{\tiny{B}}}=1.2$. Middle panel: $k^{\textrm{\tiny{T}}}_{\textrm{\tiny{R}}}=k^{\textrm{\tiny{C}}}_{\textrm{\tiny{B}}}=1.2$. Right panel: $k^{\textrm{\tiny{T}}}_{\textrm{\tiny{B}}}=k^{\textrm{\tiny{C}}}_{\textrm{\tiny{B}}}=1.2$.}
\label{fig:scan_case2_zeros}
\end{center}
\end{figure*}

The trajectories of the final region inside the white contours, centred at the point $k^{\textrm{\tiny{T}}}_{\textrm{\tiny{B}}} = k^{\textrm{\tiny{T}}}_{\textrm{\tiny{R}}} = k^{\textrm{\tiny{C}}}_{\textrm{\tiny{R}}} = k^{\textrm{\tiny{C}}}_{\textrm{\tiny{B}}} = 1.2$, the similarity of the parameter values means the trajectories for the engagements between Red and Blue display a complex back-and-forth dynamic as the population of each force increases and decreases a number of times.  Figure \ref{fig:scan_case2_zeros} displays the number of periods experienced by Red's trajectory as a function of input parameters, with the middle panel corresponding to the parameter inputs used in Figure \ref{fig:fig_case2_scan}. Focusing on the middle panel in Figure \ref{fig:scan_case2_zeros}, the region bounded by the white contours exhibits significant oscillatory behaviour --- centred at $k^{\textrm{\tiny{T}}}_{\textrm{\tiny{B}}} = k^{\textrm{\tiny{T}}}_{\textrm{\tiny{R}}} = k^{\textrm{\tiny{C}}}_{\textrm{\tiny{R}}} = k^{\textrm{\tiny{C}}}_{\textrm{\tiny{B}}} = 1.2$. Mirroring the top two panels of Figure \ref{fig:fig_case2}, the oscillations of the Red and Blue trajectories are centred around each force's carrying capacity value. There is a stark contrast when comparing the final total Green population(s) in this region with the regions above and below the white contours. The oscillations generally serve to prolong the engagement, as seen by the difference between the time-scales of the top and bottom panels of Figure \ref{fig:fig_case2}. Thus the final total Green populations in the oscillatory region is notably much less than the surrounding regions as the prolonged time period means that a greater proportion of the Green population contributes to the conflict.

\subsection{Contributing case: analytic considerations}
\label{SEC4-2}

By applying approximations in this instance of the model it is possible to gain greater understanding of the oscillatory behaviour displayed by the Red and Blue trajectories. As with the previous case in section \ref{modIanalyt} where the near-peer scenario was examined in depth, to allow analytical amenability we assign $k^{\textrm{\tiny{T}}}_{\textrm{\tiny{B}}}=k^{\textrm{\tiny{T}}}_{\textrm{\tiny{R}}}$ and $k^{\textrm{\tiny{C}}}_{\textrm{\tiny{R}}} = k^{\textrm{\tiny{C}}}_{\textrm{\tiny{B}}}$ in Eq.(\ref{Greenpop2}), and consider the case $\gamma = g$ and $B=R$. Doing so to we obtain the reduced system:
\begin{eqnarray}
\dot{B} = - B + k^{\textrm{\tiny{T}}}_{\textrm{\tiny{B}}} g, \;\;
\dot{g} = g \Gamma B \left(1 - \frac{B}{k^{\textrm{\tiny{C}}}_{\textrm{\tiny{R}}}}\right) - k^{\textrm{\tiny{T}}}_{\textrm{\tiny{B}}} g,\;\;
\dot{\Gamma} = - 2 g \Gamma B \left(1 - \frac{B}{k^{\textrm{\tiny{C}}}_{\textrm{\tiny{R}}}} \right).
\label{simplify1}
\end{eqnarray}
Further assuming that $B$ exhibits small oscillations around the value $k^{\textrm{\tiny{C}}}_{\textrm{\tiny{R}}}$, behaviour that was exhibited in the top two panels of Figure \ref{fig:fig_case2}. Hence,
\begin{equation}
    B(t) = k^{\textrm{\tiny{C}}}_{\textrm{\tiny{R}}} + \epsilon(t),
    \label{small1}
\end{equation}
where we assume that the time dependent variable in Eq.~(\ref{small1}) is small, \textit{i.e.} $\epsilon^2 \approx 0$. With this assumption, the system in Eq.~(\ref{simplify1}) becomes,
\begin{equation}
    \dot{\epsilon} = - k^{\textrm{\tiny{C}}}_{\textrm{\tiny{R}}} - \epsilon + k^{\textrm{\tiny{T}}}_{\textrm{\tiny{B}}} g, \;\;
    \dot{g} = - g \Gamma \epsilon - k^{\textrm{\tiny{T}}}_{\textrm{\tiny{B}}} g, \;\;
    \dot{\Gamma} = 2 g \Gamma \epsilon,
    \label{simplify2}
\end{equation}
where any quadratic terms in $\epsilon$ are approximately zero. The multiplicative $\epsilon$ terms on the right hand side of Eq.~(\ref{simplify2}) highlight their role in determining the direction of flow between $g$ and $\Gamma$. That is, if $\epsilon < 0$, Blue can maintain support from Green, and population will flow from $\Gamma$ to $g$. Contrastingly, if $\epsilon >0$, Blue cannot maintain support and contributions from Green, and population will flow from $g$ to $\Gamma$. 

Summing all the terms in Eq.~(\ref{simplify2}), and assuming that the initial condition for $\epsilon$ is zero (\textit{i.e.} $\epsilon_0=0$), we obtain
\begin{eqnarray}
\begin{split}
&   \dot{\epsilon}+ \dot{g}+ \frac{1}{2}\dot{\Gamma} =- k^{\textrm{\tiny{C}}}_{\textrm{\tiny{R}}} -\epsilon, \\
\Rightarrow\;\; & \epsilon+ g+\frac{1}{2}\Gamma = g_0 + \frac{1}{2}\Gamma_0 - k^{\textrm{\tiny{C}}}_{\textrm{\tiny{R}}} t - \int^t_0 d\tau \epsilon(\tau).
\label{suggestive}
\end{split}
\end{eqnarray}
In the dynamical region of interest of oscillatory $\epsilon$, we can assume the term $\int^t_0 d\tau \epsilon(\tau)$ contributes very little to Eq.(\ref{suggestive}) as the trajectory of $\epsilon$ will not have strayed significantly negatively away from zero. Thus at $t=t_f$, where
\begin{equation}
t_f= \frac{1}{k^{\textrm{\tiny{C}}}_{\textrm{\tiny{R}}}}\left( g_0 + \frac{\Gamma_0}{2}\right),
\label{final_timescale}
\end{equation}
the right hand side of the final expression of Eq.~(\ref{suggestive}) approximately equals zero, giving a natural time-scale, past which the oscillatory dynamics will cease, as either both $g$ and $\Gamma$ are zero, or $\epsilon$ is significantly below zero.

Dividing $\dot{g}$ by $\dot{\Gamma}$ in Eq.~(\ref{simplify2}) reveals
\begin{eqnarray}
\begin{split}
2 \frac{dg}{d \Gamma} = -1 - \frac{k^{\textrm{\tiny{T}}}_{\textrm{\tiny{B}}}}{\Gamma \epsilon},\\
\Rightarrow - 2 \epsilon (d \epsilon + k^{\textrm{\tiny{C}}}_{\textrm{\tiny{R}}} dt + \epsilon dt ) = - k^{\textrm{\tiny{T}}}_{\textrm{\tiny{B}}} \frac{\ d\Gamma}{\Gamma}, \\
\Rightarrow \epsilon^2 + 2 \int^t_0 d\tau \left[k^{\textrm{\tiny{C}}}_{\textrm{\tiny{R}}} \epsilon(\tau) + \epsilon^2(\tau) \right] = k^{\textrm{\tiny{T}}}_{\textrm{\tiny{B}}} \ln \frac{\Gamma}{\Gamma_0},\\
\Rightarrow \Gamma = \Gamma_0 \exp \left\{ \frac{\epsilon^2 + 2 \int^t_0 d\tau \left[k^{\textrm{\tiny{C}}}_{\textrm{\tiny{R}}} \epsilon(\tau) + \epsilon^2(\tau) \right]}{k^{\textrm{\tiny{T}}}_{\textrm{\tiny{B}}}}  \right\}
\end{split}
\label{workings1}
\end{eqnarray}
where Eq.~(\ref{suggestive}) has been applied to proceed from the first line of Eq.~(\ref{workings1}) to the second. Additionally, the expression $\int^t_0 d\tau \epsilon(\tau)$, and the terms nonlinear in $\epsilon$ in Eq.(\ref{workings1}) have been kept, with a plan to deal with them appropriately at the conclusion of this section. Eq.~(\ref{workings1}) gives the trajectory of $\Gamma$ entirely in terms of $\epsilon$. In order to utilise this expression, we re-cast the equation for $g$ in Eq.~(\ref{simplify2}) as
\begin{equation}
   e^{-k^{\textrm{\tiny{T}}}_{\textrm{\tiny{B}}} t}\frac{d}{dt}g e^{k^{\textrm{\tiny{T}}}_{\textrm{\tiny{B}}} t} = \dot{g} + k^{\textrm{\tiny{T}}}_{\textrm{\tiny{B}}} g = - \frac{1}{2}\dot{\Gamma}.
\label{workings2}
\end{equation}
Performing a similar operation to the equation for $\epsilon$ in Eq.~(\ref{simplify2}) --- \textit{i.e.} multiplying the entire expression by $e^{k^{\textrm{\tiny{T}}}_{\textrm{\tiny{B}}} t}$, and then differentiating with respect to time --- we obtain
\begin{eqnarray}
\begin{split}
    \ddot{\epsilon} + (k^{\textrm{\tiny{T}}}_{\textrm{\tiny{B}}} +1) \dot{\epsilon} + k^{\textrm{\tiny{T}}}_{\textrm{\tiny{B}}} \epsilon = -\frac{k^{\textrm{\tiny{T}}}_{\textrm{\tiny{B}}}}{2} \dot{\Gamma} - k^{\textrm{\tiny{T}}}_{\textrm{\tiny{B}}} k^{\textrm{\tiny{C}}}_{\textrm{\tiny{R}}}\\
    \Rightarrow \dot{\epsilon} + (k^{\textrm{\tiny{T}}}_{\textrm{\tiny{B}}}+1) \epsilon + k^{\textrm{\tiny{T}}}_{\textrm{\tiny{B}}} \int^t_0 d \tau \epsilon(\tau) + \frac{k^{\textrm{\tiny{T}}}_{\textrm{\tiny{B}}}}{2} \Gamma 
    = k^{\textrm{\tiny{T}}}_{\textrm{\tiny{B}}} \left( g_0 + \frac{\Gamma_0}{2} -k^{\textrm{\tiny{C}}}_{\textrm{\tiny{R}}} t \right) - k^{\textrm{\tiny{C}}}_{\textrm{\tiny{R}}},
\end{split}
\label{workings3}
\end{eqnarray}
with initial condition $\dot{\epsilon}_0 = - k^{\textrm{\tiny{C}}}_{\textrm{\tiny{R}}} + k^{\textrm{\tiny{T}}}_{\textrm{\tiny{B}}} g_0$. As the expression for $\Gamma$ in terms of $\epsilon$ in Eq.~(\ref{workings1}) contains an integral of $\epsilon$ with respect to time, the expression in Eq.~(\ref{workings3}) gives an integro-differential equation for the trajectory of $\epsilon$. To change all but the argument of the exponential of Eq.~(\ref{workings3}) into an ODE, we redefine the time dependent variables as follows,
\begin{equation}
    \alpha = \int^t_0 d \tau \epsilon(\tau), \;\; \dot{\alpha} = \epsilon, \;\; \ddot{\alpha} = \dot{\epsilon},
\end{equation}
with initial conditions $\alpha_0 = \dot{\alpha}_0 = 0$, $\ddot{\alpha}_0 = - k^{\textrm{\tiny{C}}}_{\textrm{\tiny{R}}} + k^{\textrm{\tiny{T}}}_{\textrm{\tiny{B}}} g_0.$ Doing so, Eq.~(\ref{workings3}) becomes:
\begin{equation}
    \ddot{\alpha}  + \frac{k^{\textrm{\tiny{T}}}_{\textrm{\tiny{B}}} \Gamma_0}{2} \exp \left\{  \frac{ \dot{\alpha}^2 + 2\left[k^{\textrm{\tiny{C}}}_{\textrm{\tiny{R}}} \alpha + \int^t_0 d \tau \dot{\alpha}^2(\tau) \right]}{k^{\textrm{\tiny{T}}}_{\textrm{\tiny{B}}}}\right\}
   + (k^{\textrm{\tiny{T}}}_{\textrm{\tiny{B}}} +1) \dot{\alpha} + k^{\textrm{\tiny{T}}}_{\textrm{\tiny{B}}} \alpha = k^{\textrm{\tiny{T}}}_{\textrm{\tiny{B}}}\left( g_0 + \frac{\Gamma_0}{2} - k^{\textrm{\tiny{C}}}_{\textrm{\tiny{R}}} t  \right)-k^{\textrm{\tiny{C}}}_{\textrm{\tiny{R}}}.
    \label{alphtraj}
\end{equation}
The expression for the trajectory of $\epsilon$ --- \textit{i.e.} $\dot{\alpha}$ --- in Eq.~(\ref{alphtraj}) produces exactly the trajectories for $\epsilon$ numerically generated from the system in Eq.~(\ref{simplify2}). Nevertheless, being nonlinear itself Eq.~(\ref{alphtraj}) still requires numerical integration.

\begin{figure*}[!htb]
\begin{center}
\includegraphics[width=14.3cm]{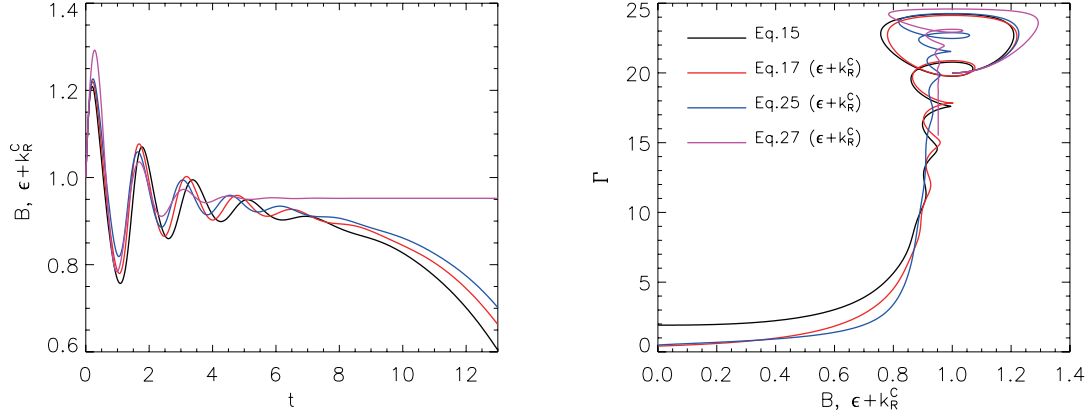}
\caption{Solutions for $B$ and $\epsilon$ (left panel), and phase portraits (right panel) for the systems Eqs.~\ref{simplify1}, \ref{simplify2}, and Eqs.~\ref{alphtraj15} and \ref{alphatraj2}, are shown by black, red, blue and magenta curves, respectively. The solution is produced for $B_0=B(t=0)=1$, $g_0=g(t=0)=3$, $\Gamma_0=\Gamma(t=0)=20$, and $k^{\textrm{\tiny{T}}}_{\textrm{\tiny{B}}}=k^{\textrm{\tiny{C}}}_{\textrm{\tiny{R}}}=1$. The solutions for $\epsilon$ are shifted by $k^{\textrm{\tiny{C}}}_{\textrm{\tiny{R}}}$ for comparison (see Eq.~\ref{small1}).}
\label{fig:spiral}
\end{center}
\end{figure*}

In order to gain analytical insight, we first approximate Eq.~(\ref{alphtraj}) into an ODE by ignoring the integral term in the exponent to obtain,
\begin{equation}
    \ddot{\alpha} + (k^{\textrm{\tiny{T}}}_{\textrm{\tiny{B}}} +1) \dot{\alpha}+ k^{\textrm{\tiny{T}}}_{\textrm{\tiny{B}}} \alpha + \frac{k^{\textrm{\tiny{T}}}_{\textrm{\tiny{B}}} \Gamma_0}{2} \exp \left( \frac{ \dot{\alpha}^2 + 2 k^{\textrm{\tiny{C}}}_{\textrm{\tiny{R}}} \alpha }{k^{\textrm{\tiny{T}}}_{\textrm{\tiny{B}}}}\right) = k^{\textrm{\tiny{T}}}_{\textrm{\tiny{B}}}\left( g_0 + \frac{\Gamma_0}{2} - k^{\textrm{\tiny{C}}}_{\textrm{\tiny{R}}} t  \right)- k^{\textrm{\tiny{C}}}_{\textrm{\tiny{R}}}.
    \label{alphtraj15}
\end{equation}
which is of the form of a \textit{forced autonomous equation} --- refer to section 2.2 of \citep{Polyanin03} for more details. Next we assume that the quadratic term of $\dot{\alpha}$ in exponential of Eq.(\ref{alphtraj15}) equals zero (\textit{i.e.} $\dot{\alpha}^2 \equiv \epsilon^2 \approx 0$). Additionally, we linearise the remaining term in the exponential via,
\begin{equation}
    \exp \left( \frac{2 k^{\textrm{\tiny{C}}}_{\textrm{\tiny{R}}}}{k^{\textrm{\tiny{T}}}_{\textrm{\tiny{B}}}}\alpha \right) \approx 1 + \frac{2 k^{\textrm{\tiny{C}}}_{\textrm{\tiny{R}}}}{k^{\textrm{\tiny{T}}}_{\textrm{\tiny{B}}}} \alpha.
\end{equation}
Doing so, Eq.~(\ref{alphtraj15}) becomes,
\begin{eqnarray}
\begin{split}
    \ddot{\alpha} + (k^{\textrm{\tiny{T}}}_{\textrm{\tiny{B}}}+1) \dot{\alpha} + (k^{\textrm{\tiny{T}}}_{\textrm{\tiny{B}}} + k^{\textrm{\tiny{C}}}_{\textrm{\tiny{R}}} \Gamma_0) \alpha = k^{\textrm{\tiny{T}}}_{\textrm{\tiny{B}}} (g_0 - k^{\textrm{\tiny{C}}}_{\textrm{\tiny{R}}} t) - k^{\textrm{\tiny{C}}}_{\textrm{\tiny{R}}}.
\end{split}
\label{alphatraj2}
\end{eqnarray}
In the left and right panels of Figure \ref{fig:spiral} we plot time-dependent trajectories for $B$ and $\epsilon$, and the $\Gamma$ vs $B$ and $\epsilon$ phase-space trajectories for the system Eq.~(\ref{simplify1}), and its various approximations given by Eqs.~(\ref{simplify2}), (\ref{alphtraj15}) and (\ref{alphatraj2}). Importantly, the solutions for $\epsilon$ are shifted by $k^{\textrm{\tiny{C}}}_{\textrm{\tiny{R}}}$ for comparison as per Eq.~(\ref{small1}). On the left panel we see that all four time-dependent trajectories initially produce very similar profiles. Nevertheless (and perhaps unsurprisingly) the trajectory corresponding to linearised Eq.(\ref{alphatraj2}) deviates significantly from the other three profiles, unable to reproduce the nonlinearity clearly being demonstrated for the entirety of the time-scale Eq.(\ref{final_timescale}).

Thus, although the linearity of Eq.(\ref{alphatraj2}) clearly limits its predictive ability over the entire time-scale, it is nonetheless immediately utilisable as the general solution is a sum of a linear combinations of exponents and a linear function. Importantly, oscillatory solutions exist if the discriminant of the characteristic equation is negative, \textit{i.e.}
\begin{equation}
k^{\textrm{\tiny{T}}}_{\textrm{\tiny{B}}} < 1+ 2 \sqrt{k^{\textrm{\tiny{C}}}_{\textrm{\tiny{R}}} \Gamma_0},    
\label{alpha_oscil}    
\end{equation}
and the oscillation period is then given as
\begin{equation}
T=\frac{2\pi}{\omega} = \frac{4\pi}{\sqrt{|( k^{\textrm{\tiny{T}}}_{\textrm{\tiny{B}}}-1)^2-4 k^{\textrm{\tiny{C}}}_{\textrm{\tiny{R}}}\Gamma_0|}}.
\label{alpha_osc_freq}
\end{equation}
Therefore, the number of oscillation periods in the solution can be approximated as the ratio of the timescale, given in Eq.~(\ref{final_timescale}), to the oscillation period as:
\begin{equation}
    N \simeq \frac{t_f}{T} = \frac{\left( 2g_0 + \Gamma_0\right)\sqrt{|(k^{\textrm{\tiny{T}}}_{\textrm{\tiny{B}}}-1)^2-4 k^{\textrm{\tiny{C}}}_{\textrm{\tiny{R}}}\Gamma_0|}}{8\pi k^{\textrm{\tiny{C}}}_{\textrm{\tiny{R}}}}.
\label{number_osc}
\end{equation}
Hence the combination of Eq.(\ref{number_osc}) and (\ref{alpha_osc_freq}) enables an approximate understanding of the effect that changes in input parameters and initial conditions have on the applicable regions, and number of oscillations occurring within the the time-scale. 

\begin{figure}
\begin{center}
\includegraphics[width=7.0cm]{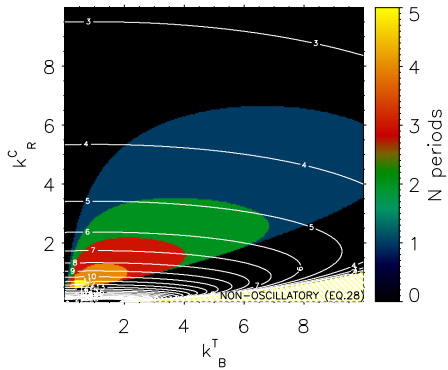}
\caption{Number of periods calculated from Eq.~\ref{number_osc} (white numbered contours) in comparison to the number of periods in the solution of the system Eq.~\ref{simplify1} at final time, given by the coloured regions ranging from black ($0$ periods) to yellow ($5$ periods). The initial conditions are as in Fig.~\ref{fig:spiral}.}
\label{fig:notagree}
\end{center}
\end{figure}

In Figure \ref{fig:notagree} we compare the oscillatory regions, and number of periods, generated from both the original system in Eq.(\ref{simplify1}) --- represented by the coloured contours --- and the fully linearised approximations given by equations~(\ref{alpha_oscil}) and (\ref{number_osc}) --- represented by the white numbered contours. The contours generated by Eq.(\ref{number_osc}) do display similar macroscopic properties exhibited by the system in Eq.(\ref{simplify1}). Notably, that an increase in oscillations is approximately inversely proportional to the carrying capacity value $k^{\textrm{\tiny{C}}}_{\textrm{\tiny{R}}}$. Additionally, we also see that as the transfer rate $k^{\textrm{\tiny{T}}}_{\textrm{\tiny{B}}}$ decreases, this generally leads to an increase in oscillation count --- though from the contours generated by Eq.(\ref{simplify1}), we see that a \textit{further} decrease in $k^{\textrm{\tiny{T}}}_{\textrm{\tiny{B}}}$ leads to the oscillation count decreasing, something not reflected in the contours of the linear system.


\section{Discussion and future work}
\label{CONC}
In this work we have extended the classical adversarial Blue vs Red model of combat attrition to include a third \textit{neutral} population who are given the opportunity to become \textit{supporters} and add to the combat effectiveness of Red and Blue, or \textit{contributors} who seek to actively join Red and Blue in their engagement. Blue or Red victory in the supporter model was most sensitive to changes of the carrying capacity. Notably, if the carrying capacity value was less than the corresponding force's initial population, the likely outcome was defeat for that force in the absence of a change of strategy.
We can see such sensitivity to carrying
capacity in the aftermath of the 2003 Iraq war. The failure to prevent
looting in the immediate days, followed by 
the later de-Ba'athification triggering
the insurgency \citep{Hosmer07}
may be interpreted 
as a collapse in carrying capacity
for US forces within its potential
supporting population. The subsequent 2007 Surge under the
`Petraeus doctrine' \citep{Petraeus2009} represented both a numerical increase in forces to defeat insurgents 
{\it and} the provision of security
to local populations to enhance
internal governance, effectively to increase the carrying capacity in potential
supportive communities.


The contributor model revealed regions of clear Blue and Red victories which surrounded a highly volatile region where outcome prediction were largely impossible. This volatility was due to both Blue and Red, and their corresponding Green supporters, being relatively evenly matched, leading to a long protracted conflict. During this conflict Blue and Red population values oscillated about their carrying capacity values, fueled by the continual contribution of Green population into both Blue and Red combatants, until finally being depleted. This scenario is reminiscent of the early stages of the 2014 conflict in Eastern Ukraine, where near-peer Ukrainian forces and pro-Russian separatists received an influx of contributors from surrounding regions \citep{Kofman17}. 

Under sufficiently symmetrised assumptions, corresponding
to a worst-case stalemate scenario, for the supporter variant we related the dynamic variables, initial conditions and the input parameters. This relationship offered a computationally inexpensive means to determine if the engagement between Red and Blue would end due to Green rescinding support completely in a near-peer setting. For the contributor model, under similar assumptions which sufficiently symmetrised defining equations, we were able to analyse the protracted oscillatory regime where a near-peer Red and Blue fluctuated around the carrying capacity value. Notably, we were able to derive a closed form expression of the expected time-scale of the oscillations, and a non-linear ODE describing the oscillations, both of which performed quite well when compared numerically with the original system. 

As the original intent of this work was to devise a relatively simple model which offered a means to explore the mutual influence between combatants and non-combatants in warfare, we shall devote some effort in interpreting these model behaviours through this lens; firstly by discussing possible meanings and implications for each force's carrying capacity value. Carrying capacity was the means we introduced each Green sub-population's growth and decline of support. If either Blue or Red grew too big, their potential Green followers retract support. This is to be interpreted as the Green population perceiving either force as acting as a military occupation, and the negative consequences which stem from that. Consequences which may be as devastating as \textit{collateral damage}, or as the \textit{negative perception} held by a local population that either Red or Blue are exerting too much influence on local decision-making. Hence, a higher carrying capacity value for either Red or Blue equates to that force being more disciplined so as not to cause collateral damage, and willing to work \textit{with} a local population, allowing for self-determination and not denying their ability for meaningful decision-making. Such a force would possess superior training, both military and cultural, and be steeped in doctrine which does not shy away from distributed and/or collaborative decision-making.

Bringing this back to behaviours witnessed by the model, we see that if it is Red and/or Blue's intent to ensure a decisive victory, then figures \ref{fig:fig_case1_scan} and \ref{fig:fig_case2_scan} clearly show that the method to ensure this through possessing a superior carrying capacity (training, doctrine, decision-making, trust, \textit{etc.}) for either supporter or contributor model variant, as this parameter is the most sensitive to change in model outcomes once approximate parity is reached. In the supporter model, with ensured conservation of Green, it would be sensible for either Blue or Red to aim for advantageous regions if there is approximate parity in carrying capacities --- equivalent to the band of Red victory surrounded by Blue victories on either side in figure \ref{fig:fig_case1_scan}. If however there is possibility for Green involvement in the engagement, similar to the contributor model variant, then it is paramount to avoid approximate parity in carrying capacity values as this leads to protracted engagements with potentially disastrous consequences on the Green population as witnessed by the depleted population values in figure \ref{fig:fig_case2_scan}.

For future applications, an immediate generalisation would include a mixed variant of our model, with one force being supported and one being contributed to, similar to Deitchman's (\citeyear{Deitchman62}) asymmetric guerrilla warfare model for the original Lanchester equations. It is our intention nevertheless to couple this model with a system of networked decision-makers via the Kuramoto-Sakaguchi model considered in \citep{Zuparic?, Kalloniatis20}. Ideally such a coupled model would possess the complexities exhibited in the model considered in this work, with the advantage of having the ability to more effectively unpack the interpretation of the carrying capacity term into its implied combat and decision-making meanings.

\section*{Acknowledgements}
This research was a collaboration between the Commonwealth of Australia (represented by the Defence Science and Technology Group) and Deakin University through a Defence Science Partnerships agreement, under the auspices of the Modelling Complex Warfighting
initiative.



\end{document}